\documentclass[11pt,reqno]{amsart}
\usepackage{amsmath,amssymb,amsthm,mathtools}
\usepackage[T1]{fontenc}
\usepackage{lmodern,microtype}
\usepackage{enumitem}
\usepackage[hidelinks]{hyperref}
\usepackage{xurl}
\usepackage{xcolor}
\numberwithin{equation}{section}
\newtheorem{theorem}{Theorem}[section]
\newtheorem{proposition}[theorem]{Proposition}
\newtheorem{lemma}[theorem]{Lemma}
\theoremstyle{definition}
\newtheorem{definition}[theorem]{Definition}
\theoremstyle{plain}
\newtheorem{corollary}[theorem]{Corollary}
\newtheorem{conjecture}[theorem]{Conjecture}
\newcommand{\M}{\mathcal M}
\newcommand{\W}{\mathcal W}
\newcommand{\Ssh}{\mathsf S}
\newcommand{\tr}{\operatorname{tr}}
\newcommand{\Tr}{\operatorname{Tr}}
\newcommand{\spec}{\operatorname{spec}}
\newcommand{\eps}{\varepsilon}
\newcommand{\C}{\mathbb C}
\newcommand{\R}{\mathbb R}
\newcommand{\N}{\mathbb N}
\title[Fixed-occurrence shuffle models]{Fixed-occurrence cases of the shuffle-word-embedding conjecture}
\author{Trung Hoa Dinh}
 \thanks{Laboratory for Applied and Industrial Mathematics, Institute for Computational Science and Artificial Intelligence, Van Lang University, Ho Chi Minh City, Vietnam; Faculty of Basic Sciences, Van Lang University, Ho Chi Minh City, Vietnam. Email: hoa.dinhtrung@vlu.edu.vn.}
 \thanks{Department of Mathematics and Statistics, Troy University, USA}
\begin{document}
\begin{abstract}
\emergencystretch=1.5em
For self-adjoint contractions in a finite tracial von Neumann algebra, binary shuffle traces with $m=1,2,4$ occurrences of the second variable admit simultaneous finite-dimensional approximation by commuting diagonal contractions. For every odd $m$, a commuting diagonal self-adjoint pair reproduces any prescribed finite set of shuffle traces exactly. When the first operator has finite spectrum and $m\geq3$ is odd, a weighted symmetric cycle reproduces all degrees simultaneously. We prove fourth-order trace-derivative positivity in finite tracial von Neumann algebras and obtain a positive measure representing the four-occurrence coefficients. For each even $m$, we characterize the contraction model property by positivity of an associated trace-derivative functional.
\end{abstract}
\maketitle

\section{Introduction}
The approximation of tracial moments by finite matrices is a basic problem
in operator algebras. The negative solution to the Connes embedding problem
shows that such approximation is not always possible for arbitrary
noncommutative words. It remains natural, however, to ask whether particular
families of tracial expressions admit finite-dimensional models. Pascoe and
Woerdeman formulated the shuffle-word-embedding conjecture in this direction:
instead of retaining the traces of individual words, one considers their
symmetrized sums \cite{PascoeWoerdeman}. The conjecture is closely related
to the Bessis--Moussa--Villani (BMV) polynomials, and asks for simultaneous
matrix approximation of the corresponding traces for self-adjoint
contractions in a type $\mathrm{II}_1$ factor.

The BMV conjecture was proposed by Bessis, Moussa and Villani
\cite{BessisMoussaVillani}, reformulated in terms of shuffle-polynomial
coefficients by Lieb and Seiringer \cite{LiebSeiringer}, and proved by
Stahl \cite{Stahl}. Collins, Dykema and Torres-Ayala studied cyclic
sum-of-squares methods for these polynomials \cite{CollinsDykemaTorres}.
A related line of work concerns the effect of spectral pinching on
symmetrized words: a pinching-corrected BMV inequality was proposed in
\cite{DinhPinching} and established for all symmetrized word averages in
\cite{DinhLeVuongBMV}. Whereas those results compare traces before and
after pinching, the present paper addresses their finite-dimensional
realization.

Let $(\M,\tau)$ be a finite von Neumann algebra with a normal faithful
tracial state satisfying $\tau(1)=1$. For self-adjoint $A,B\in\M$, define
\begin{equation*}
 \Ssh_{n,m}(A,B)=\sum_{W\in\W_{n,m}}W(A,B),\qquad n,m\in\N_0,
\end{equation*}
where $\W_{n,m}$ denotes the set of words containing $n$ occurrences of
$A$ and $m$ occurrences of $B$. Equivalently, $\Ssh_{n,m}(A,B)$ is the
coefficient of $s^nt^m$ in $(sA+tB)^{n+m}$. For matrices, write
$\tr_d(X)=d^{-1}\Tr(X)$. We fix the occurrence number $m$ and ask whether,
for any finite set $N\subset\N_0$, the traces
\[
 \tau\bigl(\Ssh_{n,m}(A,B)\bigr),\qquad n\in N,
\]
can be approximated simultaneously by one finite-matrix pair.

\begin{definition}\label{lem:model-def}
Fix $m\geq1$. The \emph{contraction model property} holds at occurrence number $m$ if, for every finite tracial von Neumann algebra $(\M,\tau)$, self-adjoint contractions $A,B\in\M$, finite set $N\subset\N_0$, and $\eps>0$, there are $d\geq1$ and self-adjoint contractions $X,Y\in M_d(\C)$ such that
\begin{equation}\label{eq:model}
 \left|\tau\bigl(\Ssh_{n,m}(A,B)\bigr)-\tr_d\bigl(\Ssh_{n,m}(X,Y)\bigr)\right|<\eps
 \quad(n\in N).
\end{equation}
The \emph{self-adjoint model property} requires $X=X^*$ and $Y=Y^*$ in \eqref{eq:model}, without the contraction bounds.
\end{definition}

Our main result establishes several fixed-occurrence cases of the
shuffle-word-embedding conjecture. In each case the approximating
matrices may be chosen to commute, although the norm requirement on the
second matrix differs between even and odd occurrences.

\begin{theorem}\label{thm:main}
Let $(\M,\tau)$ be any finite tracial von Neumann algebra and let
$A,B\in\M$ be self-adjoint contractions.
\begin{enumerate}[label=\textup{(\roman*)}]
\item At $m=1,2,4$, for every finite $N\subset\N_0$ and every
$\eps>0$, there are commuting diagonal self-adjoint contractions $X,Y$
such that
\[
\bigl|\tau(\Ssh_{n,m}(A,B))-
\tr_d(\Ssh_{n,m}(X,Y))\bigr|<\eps\quad(n\in N).
\]
\item For every odd $m\geq3$ and every finite $N\subset\N_0$,
there are commuting diagonal self-adjoint matrices $X,Y$, with
$\|X\|\leq1$, such that
\[
\tau(\Ssh_{n,m}(A,B))=
\tr_d(\Ssh_{n,m}(X,Y))\quad(n\in N).
\]
\end{enumerate}
\end{theorem}

The four-occurrence case admits a stronger description. We establish
fourth-order trace-derivative positivity for self-adjoint elements of any
finite tracial von Neumann algebra: if $f$ is a real polynomial with
$f^{(4)}\geq0$ on an interval containing $\spec(A)$, then
\[
 \left.\frac{d^4}{dt^4}\tau f(A+tB)\right|_{t=0}\geq0.
\]
The corresponding matrix inequality was proved by Hein\"avaara
\cite[Theorem~1.0.4]{HeinavaaraPlanes}; his later work extends
matrix trace-derivative positivity to every even order
\cite[Corollary~1.2]{HeinavaaraTJSM}. Trace-derivative measures are
also studied in higher-order spectral shift theory
\cite{DykemaSkripka,PradhanSkripka}. Our proof in finite tracial von
Neumann algebras leads to a positive Borel measure $\nu_{A,B}^{(4)}$
on $[-1,1]$ satisfying
\[
 \tau\bigl(\Ssh_{n,4}(A,B)\bigr)
 =\binom{n+4}{4}\int_{-1}^1 x^n\,d\nu_{A,B}^{(4)}(x),
 \qquad n\geq0,
\]
with $\nu_{A,B}^{(4)}([-1,1])=\tau(B^4)$. Thus the normalized
four-occurrence coefficients form a positive moment sequence; in
particular, their finite Hankel matrices are positive semidefinite.
Approximating the measure by rational atomic measures yields the
commuting diagonal contraction models in the theorem.

The measure representation also relates the four-occurrence result to
the analytic functions considered by Pascoe and Woerdeman. Setting
$F_{A,B}(s,t)=\exp\bigl(\tau\log(1-sA-tB)\bigr)$, we obtain
\[
 [t^4]\log F_{A,B}(s,t)
 =-\frac14\int_{-1}^1\frac{d\nu_{A,B}^{(4)}(x)}{(1-sx)^4},
 \qquad |s|<1.
\]
This identity is established in Corollary~\ref{cor:four-log}.
In particular, the fixed-occurrence approximations hold without
assuming the Connes embedding property for $\M$; they do not assert
approximation of arbitrary noncommutative moments.

At odd occurrence numbers, the exact diagonal models follow from
finite interpolation, while a weighted symmetric cycle realizes all
degrees when $A$ has finite spectrum. For even $m$, define
\begin{equation}\label{eq:Lambda-intro}
 \Lambda_{A,B}^{(m)}(q)=\left.\frac{d^m}{dt^m}\tau f(A+tB)\right|_{t=0},
 \qquad f^{(m)}=q,
\end{equation}
for real polynomials $q,f$. This does not depend on the choice of $f$.
We show that the contraction model property at even $m$ is equivalent
to $\Lambda_{A,B}^{(m)}(q)\geq0$ for every real polynomial
$q\geq0$ on $[-1,1]$. The first remaining even case is $m=6$; we
exhibit an indefinite kernel arising from a natural decomposition of
the six-occurrence terms.

{\bf Use of Artificial Intelligence.}  The author used ChatGPT to explore connections between the pinching-corrected BMV conjecture, the Connes embedding problem, and the shuffle-word-embedding conjecture of Pascoe and Woerdeman. ChatGPT also assisted with literature searches, exploratory calculations, and the investigation of possible proof strategies. The mathematical statements, arguments, proofs, and final presentation were reviewed, verified, and revised by the author, who takes full responsibility for the content of the manuscript.

The paper is organized as follows. Section~2 collects three preliminary
lemmas concerning cyclic trace identities, finite-spectrum expansions,
and atomic approximation, and then treats the cases $m=1,2$.
Section~3 proves fourth-order trace positivity and the resulting
moment and logarithmic representations. Section~4 constructs models
for odd occurrences. Section~5 gives the positivity criterion for
even occurrences, and Section~6 examines the sixth-order question.
The Appendix contains the details of the fourth-order sign calculation.

\section{Preliminaries and low-order cases}\label{sec:expansions}

\subsection{Preliminary identities}
\begin{lemma}\label{lem:cyclic}
For $m\geq1$ and $n\geq0$,
\begin{equation}\label{eq:cyclic}
 \tau\bigl(\Ssh_{n,m}(A,B)\bigr)
 =\frac{n+m}{m}\sum_{\substack{r_1,\ldots,r_m\geq0\\r_1+\cdots+r_m=n}}
 \tau\bigl(BA^{r_1}BA^{r_2}\cdots BA^{r_m}\bigr).
\end{equation}
\end{lemma}
\begin{proof}
Put $L=n+m$. Consider the set of pairs $(w,j)$ where $w$ is a word of length $L$ containing the prescribed letters and $j$ selects one of its $m$ positions occupied by $B$. Its cardinality is $m\binom{L}{m}$. Read the marked word cyclically, beginning at the selected $B$, and let $r_j$ be the number of $A$'s following the $j$th $B$. Then $r_1+\cdots+r_m=n$.

Conversely, fix such a weak composition and distinguish its first $B$. Place this distinguished cyclic word into the $L$ labeled positions of a linear word, keeping track of the distinguished position. There are $L$ choices of that position. For a periodic word, the $L$ marked positions are still distinct. Each of these $L$ choices has the same trace, by cyclicity. Summing over all compositions, and comparing with the $m$ possible markings of each unmarked word, gives $m\tau(\Ssh_{n,m})=L\sum_r\tau(BA^{r_1}\cdots BA^{r_m})$.
\end{proof}

For $x_1,\ldots,x_m\in\R$ let $h_n(x_1,\ldots,x_m)$ be the complete homogeneous polynomial,
\begin{equation*}
h_n(x_1,\ldots,x_m)=\sum_{r_1+\cdots+r_m=n}x_1^{r_1}\cdots x_m^{r_m},\qquad h_0=1.
\end{equation*}
\begin{lemma}\label{lem:finite-spectrum}
Suppose $A=\sum_{i=1}^q a_i P_i,$ 
where the $a_i$ are distinct real numbers and the nonzero $P_i$ are orthogonal projections with $\sum_iP_i=1$. Then
\begin{equation}\label{eq:finite}
 \tau\bigl(\Ssh_{n,m}(A,B)\bigr)=\frac{n+m}{m}
 \sum_{i_1,\ldots,i_m=1}^q
 h_n(a_{i_1},\ldots,a_{i_m})
 \tau(P_{i_1}BP_{i_2}B\cdots P_{i_m}B).
\end{equation}
\end{lemma}
\begin{proof}
Insert $A^{r_j}=\sum_i a_i^{r_j}P_i$ into \eqref{eq:cyclic}. For a fixed ordered tuple $(i_1,\ldots,i_m)$, summation over $r_1+\cdots+r_m=n$ gives the scalar factor $h_n(a_{i_1},\ldots,a_{i_m})$. Cyclicity of $\tau$ gives the order of factors in \eqref{eq:finite}.
\end{proof}

\begin{lemma}\label{lem:atomic}
Let $K$ be a compact subset of a Euclidean space and let $\mu$ be a positive measure on $K$ with total mass at most one. Given finitely many continuous real-valued functions $\phi_1,\ldots,\phi_s$ and $\delta>0$, there are points $z_1,\ldots,z_r\in K$, integers $k_1,\ldots,k_r\geq1$ and $d\geq\sum_jk_j$ such that
\begin{equation*}
 \left|\int_K\phi_\ell\,d\mu-\frac1d\sum_{j=1}^r k_j\phi_\ell(z_j)\right|<\delta
 \quad(1\leq\ell\leq s).
\end{equation*}
\end{lemma}
\begin{proof}
By uniform continuity, partition $K$ into finitely many Borel sets of sufficiently small diameter and replace the restriction of $\mu$ to each set by a point mass of the same weight. This changes the prescribed integrals by an arbitrarily small amount. Next choose a common denominator $d$ and approximate the finitely many nonnegative masses by nonnegative integer multiples of $1/d$, preserving the inequality that their sum is at most one. One can do this, for instance, by rounding all masses downward with a sufficiently large denominator. The total change in each integral is controlled by the sum of the rounding errors times $\max_\ell\|\phi_\ell\|_\infty$.

\end{proof}

\subsection{One occurrence}\label{sec:one}
\begin{proposition}\label{prop:one}
At $m=1$, every finite set of shuffle traces admits commuting diagonal contraction models.
\end{proposition}
\begin{proof}
Let $\mathcal A=W^*(A)$ and let $E_{\mathcal A}:\M\to\mathcal A$ be the trace-preserving conditional expectation. The defining property of $E_{\mathcal A}$ gives
\begin{equation*}
 \tau\bigl(\Ssh_{n,1}(A,B)\bigr)=(n+1)\tau(A^nB)
 =(n+1)\tau\bigl(A^nE_{\mathcal A}(B)\bigr).
\end{equation*}
The conditional expectation is positive, unital and $*$-preserving, so $E_{\mathcal A}(B)$ is again a self-adjoint contraction. It commutes with $A$. The joint spectral measure $\rho$ of $(A,E_{\mathcal A}(B))$ satisfies $\tau(A^nE_{\mathcal A}(B))=\int_{[-1,1]^2}x^ny\,d\rho(x,y)$. Apply Lemma \ref{lem:atomic} to the finitely many functions $(n+1)x^ny$, $n\in N$. Repeat the first and second coordinates of each selected atom on the diagonals of $X$ and $Y$ with the prescribed integer multiplicities, and set both entries equal to zero in the remaining coordinates. Both matrices are self-adjoint contractions. The resulting normalized traces satisfy \eqref{eq:model} for every $n\in N$.
\end{proof}

\subsection{Two occurrences}\label{sec:two}
\begin{theorem}\label{thm:two}
Let $A,B$ be self-adjoint contractions in a finite tracial von
Neumann algebra $(\M,\tau)$. There is a finite positive Borel
measure $\nu^{(2)}_{A,B}$ on $[-1,1]$ satisfying
\begin{equation}\label{eq:two-moment}
 \tau\bigl(\Ssh_{n,2}(A,B)\bigr)
 =\binom{n+2}{2}\int_{-1}^1 z^n\,d\nu^{(2)}_{A,B}(z),
 \qquad n\geq0,
\end{equation}
and $\nu^{(2)}_{A,B}([-1,1])=\tau(B^2)\leq1$.
The contraction model property at $m=2$ holds with commuting diagonal
self-adjoint contractions.
\end{theorem}
\begin{proof}
On the Hilbert space $L_2(\M,\tau)$, with
$\langle\xi,\eta\rangle=\tau(\eta^*\xi)$, let
$L_A\xi=A\xi$ and $R_A\xi=\xi A$.
These are commuting self-adjoint contractions. Their joint spectral
measure $\mu$ at the vector $B\in L_2(\M,\tau)$ is positive, supported
on $[-1,1]^2$, and satisfies
\begin{equation*}
 \int_{[-1,1]^2}x^ry^s\,d\mu(x,y)
 =\tau(BA^rBA^s),\qquad
 \mu([-1,1]^2)=\tau(B^2)\leq1.
\end{equation*}
The cyclic trace formula \eqref{eq:cyclic} gives
\begin{equation}\label{eq:two-h}
 \tau\bigl(\Ssh_{n,2}(A,B)\bigr)
 =\frac{n+2}{2}\int_{[-1,1]^2}h_n(x,y)\,d\mu(x,y).
\end{equation}
The beta-integral identity is
\begin{equation}\label{eq:two-beta}
 \int_0^1\bigl(tx+(1-t)y\bigr)^n\,dt
 =\frac1{n+1}h_n(x,y).
\end{equation}
Indeed, the coefficient of $x^ry^{n-r}$ in the integral is
$\binom nr\int_0^1t^r(1-t)^{n-r}\,dt=1/(n+1)$.
Let $\nu^{(2)}_{A,B}$ be the pushforward of
$\mu\otimes dt$ under the map $(x,y,t)\mapsto tx+(1-t)y$.
This map takes $[-1,1]^2\times[0,1]$ into $[-1,1]$, and the
pushforward has total mass $\tau(B^2)$. Combining
\eqref{eq:two-h} and \eqref{eq:two-beta} yields
\[
 \tau\bigl(\Ssh_{n,2}(A,B)\bigr)
 =\frac{(n+2)(n+1)}2\int z^n\,d\nu^{(2)}_{A,B}(z),
\]
which proves the moment formula. 

Fix a finite set $N$ and $\eps>0$. Apply Lemma \ref{lem:atomic}
to $\nu^{(2)}_{A,B}$ and the continuous functions
$\binom{n+2}{2}z^n$, $n\in N$.
Write the resulting rational atomic subprobability measure as
$\sum_j(k_j/d)\delta_{z_j}$, with $\sum_jk_j\leq d$.
Choose diagonal matrices $X,Y$ with $k_j$ copies of $(z_j,1)$
on their respective diagonals for each $j$ and with $(0,0)$ on
all remaining diagonal positions. They are commuting self-adjoint
contractions, and
\[
 \tr_d\bigl(\Ssh_{n,2}(X,Y)\bigr)
 =\binom{n+2}{2}\sum_j\frac{k_j}{d}z_j^n,
\]
which approximates \eqref{eq:two-moment} within $\eps$ for
all $n\in N$.
\end{proof}

\section{Fourth-order trace positivity}\label{sec:four}
For a $C^k$ function $g$, denote its divided difference of order $k$ by $[x_0,\ldots,x_k]_g$, with repeated nodes defined by continuity. Thus $[x,\ldots,x]_g=g^{(k)}(x)/k!$.
\begin{lemma}\label{lem:traceDD}
If $A=\sum_i\lambda_iP_i$ has finite spectrum and $f$ is a real polynomial, then
\begin{equation}\label{eq:traceDD}
 \frac1{4!}\left.\frac{d^4}{dt^4}\tau f(A+tB)\right|_{t=0}
 =\frac14\sum_{i,j,k,l}
 [\lambda_i,\lambda_j,\lambda_k,\lambda_l]_{f'}
 \tau(P_iBP_jBP_kBP_lB).
\end{equation}
\end{lemma}
\begin{proof}
Both sides depend linearly on $f$, so it suffices to take $f(x)=x^p$. When $p<4$, both sides vanish. For $p\geq4$, expanding $(A+tB)^p$ and collecting the fourth-order coefficient gives
\[
 \frac1{4!}\left.\frac{d^4}{dt^4}\tau(A+tB)^p\right|_{t=0}
 =\frac p4\sum_{r_1+\cdots+r_4=p-4}
 \tau(BA^{r_1}BA^{r_2}BA^{r_3}BA^{r_4}),
\]
by Lemma \ref{lem:cyclic}. Inserting the spectral resolution, the coefficient associated with the four indices $i,j,k,l$ is
\[
 \frac p4 h_{p-4}(\lambda_i,\lambda_j,\lambda_k,\lambda_l).
\]
The monomial divided-difference identity gives
\[
 [x_1,x_2,x_3,x_4]_{x^{p-1}}=h_{p-4}(x_1,x_2,x_3,x_4).
\]
Since $f'(x)=p x^{p-1}$, these coefficients agree with \eqref{eq:traceDD}.
\end{proof}

For $c\notin\spec(A)$, formula \eqref{eq:traceDD} applies to $g_c(x)=|x-c|^3$. In that case $g_c$ is analytic in a neighborhood of the finite spectrum, and the divided-difference trace formula remains valid: polynomial Hermite interpolation matches the values of $g_c$ and its first four derivatives at the spectral points. The formula uses only derivatives at spectral points, none of which equals $c$.

\begin{lemma}\label{lem:four-square}
Let $A=\sum_i\lambda_iP_i$ be self-adjoint with finite spectrum, with no zero eigenvalue, and let $B=B^*\in\M$. Then
\begin{equation*}
 \left.\frac{d^4}{dt^4}\tau|A+tB|^3\right|_{t=0}\geq0.
\end{equation*}
The derivative is well defined because $A+tB$ remains invertible for sufficiently small $t$.
\end{lemma}
\begin{proof}
Write $I_+=\{i:\lambda_i>0\}$ and $I_-=\{i:\lambda_i<0\}$, and set $g(x)=x|x|$. Since $(|x|^3)'=3g(x)$, equation \eqref{eq:traceDD} shows that the coefficient of $t^4$ in $\tau|A+tB|^3$ is
\begin{equation*}
 \frac34\,T_g,
 \qquad T_g=\sum_{i,j,k,l}
 [\lambda_i,\lambda_j,\lambda_k,\lambda_l]_g
 \tau(P_iBP_jBP_kBP_lB).
\end{equation*}
For $a_1,a_2,a_3<0<b_1,b_2,b_3$, the divided differences are
\begin{align}
 [a_1,a_2,a_3,b_1]_g
 &=\frac{2b_1^2}{(b_1-a_1)(b_1-a_2)(b_1-a_3)},\label{eq:DD1}\\
 [a_1,a_2,b_1,b_2]_g
 &=\frac{2(b_1+b_2)a_1a_2-2(a_1+a_2)b_1b_2}
 {(b_1-a_1)(b_1-a_2)(b_2-a_1)(b_2-a_2)},\label{eq:DD2}\\
 [a_1,b_1,b_2,b_3]_g
 &=\frac{2a_1^2}{(b_1-a_1)(b_2-a_1)(b_3-a_1)}.\label{eq:DD3}
\end{align}
The divided difference vanishes if all four nodes have the same sign, since $g$ is quadratic on each half-line. Equations \eqref{eq:DD1}--\eqref{eq:DD3} follow from the divided-difference recursion; repeated nodes follow by continuity.

Classify the nonzero summands by the cyclic sign patterns of their four spectral indices. Define, for $i\in I_+$ and $j\in I_-$,
\begin{equation}\label{eq:C}
 C_{ij}=\sum_{r\in I_-}\frac{\lambda_i}{\lambda_i-\lambda_r}
 P_iBP_rBP_j
 +\sum_{r\in I_+}\frac{\lambda_j}{\lambda_j-\lambda_r}
 P_iBP_rBP_j.
\end{equation}
For $i,j\in I_-$ and $i,j\in I_+$, respectively, define
\begin{align*}
 D^-_{ij}&=\sum_{r\in I_+}
 \frac{\lambda_r}{(\lambda_i-\lambda_r)(\lambda_j-\lambda_r)}P_iBP_rBP_j,\\
 D^+_{ij}&=\sum_{r\in I_-}
 \frac{\lambda_r}{(\lambda_i-\lambda_r)(\lambda_j-\lambda_r)}P_iBP_rBP_j.
\end{align*}
All denominators are nonzero because their indices lie in opposite sign classes. Now put
\begin{align}
 Q_{\rm adj}&=2\sum_{i\in I_+,\,j\in I_-}
 \frac{\tau(C_{ij}C_{ij}^*)}{\lambda_i-\lambda_j},\label{eq:Qadj}\\
 Q_{\rm alt}&=\sum_{i,j\in I_-}(-\lambda_i-\lambda_j)
 \tau(D^-_{ij}(D^-_{ij})^*)
 +\sum_{i,j\in I_+}(\lambda_i+\lambda_j)
 \tau(D^+_{ij}(D^+_{ij})^*).\label{eq:Qalt}
\end{align}
Every scalar coefficient in \eqref{eq:Qadj}--\eqref{eq:Qalt} is nonnegative, and $\tau(ZZ^*)\geq0$; hence $Q_{\rm adj},Q_{\rm alt}\geq0$.

The all-positive and all-negative terms vanish. Appendix~\ref{app:signs} gives the exact identities
\[
\sum_{\text{two sign changes}}(\cdots)=4Q_{\rm adj},
\qquad
\sum_{\text{alternating signs}}(\cdots)=4Q_{\rm alt},
\]
where each sum is taken over the corresponding terms of $T_g$.

Therefore,
\begin{equation}\label{eq:Gram}
 T_g=4\bigl(Q_{\rm adj}+Q_{\rm alt}\bigr),\qquad
 \frac1{4!}\left.\frac{d^4}{dt^4}\tau|A+tB|^3\right|_{t=0}
 =3\bigl(Q_{\rm adj}+Q_{\rm alt}\bigr)\geq0.
\end{equation}

\end{proof}

\begin{theorem}\label{thm:four-positive}
Let $A,B$ be self-adjoint elements of a finite tracial von Neumann algebra. If $f$ is a real polynomial and $f^{(4)}\geq0$ on a compact interval containing $\spec(A)$, then
\begin{equation*}
 \left.\frac{d^4}{dt^4}\tau f(A+tB)\right|_{t=0}\geq0.
\end{equation*}
\end{theorem}
\begin{proof}
First suppose that $A$ has finite spectrum, and let $I=[a,b]$ contain it. Taylor's formula at $a$ gives, for $x\in I$,
\begin{equation}\label{eq:Taylor4}
 f(x)=p_3(x)+\frac16\int_a^b(x-c)_+^3f^{(4)}(c)\,dc,
\end{equation}
where $p_3$ is a polynomial of degree at most three. Since
\begin{equation}\label{eq:halfabs}
 (x-c)_+^3=\frac12\bigl((x-c)^3+|x-c|^3\bigr),
\end{equation}
the integrand differs, for the purposes of a fourth-order trace derivative, from $|x-c|^3/2$ by a cubic polynomial. The finite set $\spec(A)$ has Lebesgue measure zero. For every other $c$, Lemma \ref{lem:four-square} applies to $A-c1$ and $B$. Equation \eqref{eq:traceDD} depends only on the corresponding divided differences at the spectral points. Equations \eqref{eq:Taylor4}--\eqref{eq:halfabs} allow the resulting finite sum of divided differences to be integrated in $c$, without differentiating a nonsmooth trace function at a singular spectral value. The coefficient of $t^4$ in $\tau f(A+tB)$ is the integral of the nonnegative coefficients from Lemma~\ref{lem:four-square} against $f^{(4)}(c)\,dc/12$. It is therefore nonnegative.

For general $A$, take finite-spectrum approximants $A_k=\varphi_k(A)$ by bounded real-valued step functions, choosing the step values in $\spec(A)$. Then $A_k=A_k^*$, $\|A_k-A\|\to0$, and $\spec(A_k)\subset I$. The inequality already proved applies to every $A_k$. As $f$ is a polynomial, its fourth trace derivative at zero is a finite sum of traces of words in $A_k$ and $B$. Each word is norm-continuous in $A_k$, so these derivatives converge to the derivative with $A$. The limit of nonnegative numbers is nonnegative.
\end{proof}

\begin{corollary}\label{cor:four-measure}
Let $A,B$ be self-adjoint contractions in a finite tracial von Neumann algebra. There is a positive Borel measure $\nu_{A,B}^{(4)}$ on $[-1,1]$ with
\begin{equation}\label{eq:four-measure}
 \tau\bigl(\Ssh_{n,4}(A,B)\bigr)
 =\binom{n+4}{4}\int_{-1}^1x^n\,d\nu_{A,B}^{(4)}(x)
 \quad(n\geq0),
\end{equation}
and $\nu_{A,B}^{(4)}([-1,1])=\tau(B^4)\leq1$. For every finite set $N$ and $\eps>0$, the corresponding coefficients admit commuting diagonal contraction models satisfying \eqref{eq:model}.
\end{corollary}
\begin{proof}
For a polynomial $q$, define $\Lambda_{A,B}^{(4)}(q)$ by \eqref{eq:Lambda-intro}. Theorem \ref{thm:four-positive} implies positivity whenever $q\geq0$ on $[-1,1]$. Its value at the constant polynomial $1$ is $\tau(B^4)$, because we may choose $f(x)=x^4/4!$. For any real polynomial $q$, the inequalities
$-\|q\|_\infty\leq q\leq\|q\|_\infty$ imply
\[
 |\Lambda_{A,B}^{(4)}(q)|\leq\tau(B^4)\|q\|_\infty.
\]
The functional extends continuously from polynomials to $C([-1,1])$. By the Riesz representation theorem, it is integration against a positive measure $\nu_{A,B}^{(4)}$. Its total mass is $\Lambda_{A,B}^{(4)}(1)=\tau(B^4)$.

For $n\geq0$ choose $f_n(x)=n!x^{n+4}/(n+4)!$, so $f_n^{(4)}(x)=x^n$. Expanding $f_n(A+tB)$ in $t$ gives
\[
 \Lambda_{A,B}^{(4)}(x^n)=\frac{n!\,4!}{(n+4)!}
 \tau\bigl(\Ssh_{n,4}(A,B)\bigr)
 =\binom{n+4}{4}^{-1}\tau\bigl(\Ssh_{n,4}(A,B)\bigr),
\]
which proves \eqref{eq:four-measure}.

For a finite set $N$, apply Lemma \ref{lem:atomic} to the functions $\binom{n+4}{4}x^n$, $n\in N$, obtaining the atomic measure $\sum_j(k_j/d)\delta_{x_j}$ with $\sum_jk_j\leq d$. Put $x_j$ on the diagonal of $X$ with multiplicity $k_j$, and put the corresponding entries of $Y$ equal to one. Fill the remaining diagonal positions with $X=0$ and $Y=0$. Since $X$ and $Y$ commute,
\[
 \Ssh_{n,4}(X,Y)=\binom{n+4}{4}X^nY^4.
\]
Hence $\tr_d(\Ssh_{n,4}(X,Y))=\binom{n+4}{4}\sum_j(k_j/d)x_j^n$, which differs from $\tau(\Ssh_{n,4}(A,B))$ by less than $\eps$ for each $n\in N$.
\end{proof}

\begin{corollary}\label{cor:four-log}
Let $A,B$ be self-adjoint contractions in a finite tracial von Neumann
algebra, and set
\[
 F_{A,B}(s,t)=\exp\bigl(\tau\log(1-sA-tB)\bigr),
 \qquad |s|+|t|<1.
\]
With the measure $\nu_{A,B}^{(4)}$ from Corollary~\ref{cor:four-measure},
for every $|s|<1$ we have
\[
 [t^4]\log F_{A,B}(s,t)
 =-\frac14\int_{-1}^1\frac{d\nu_{A,B}^{(4)}(x)}{(1-sx)^4},
\]
and consequently
\[
 \left.\frac{\partial^4}{\partial t^4}
 \tau\log(1-sA-tB)\right|_{t=0}
 =-6\int_{-1}^1\frac{d\nu_{A,B}^{(4)}(x)}{(1-sx)^4}.
\]
\end{corollary}
\begin{proof}
The norm-convergent power series for the logarithm gives, for
$|s|+|t|<1$,
\[
 \tau\log(1-sA-tB)
 =-\sum_{n+m\geq1}\frac{s^nt^m}{n+m}
    \tau\bigl(\Ssh_{n,m}(A,B)\bigr).
\]
Therefore, by \eqref{eq:four-measure},
\begin{align*}
 [t^4]\log F_{A,B}(s,t)
 &=-\sum_{n=0}^{\infty}\frac{s^n}{n+4}
     \tau\bigl(\Ssh_{n,4}(A,B)\bigr)\\
 &=-\frac14\int_{-1}^1
    \sum_{n=0}^{\infty}\binom{n+3}{3}(sx)^n
      \,d\nu_{A,B}^{(4)}(x)\\
 &=-\frac14\int_{-1}^1
     \frac{d\nu_{A,B}^{(4)}(x)}{(1-sx)^4}.
\end{align*}
The series converges uniformly in $x\in[-1,1]$ for $|s|<1$,
justifying the interchange of sum and integral. Multiplying the
coefficient identity by $4!=24$ gives the derivative formula.
\end{proof} 

\section{Odd occurrences: exact diagonal models and weighted cycles}\label{sec:odd}
\begin{theorem}\label{thm:odd-diagonal}
Let $m\geq1$ be odd, and let $A,B$ be self-adjoint contractions
in a finite tracial von Neumann algebra.
For every finite $N\subset\N_0$ there are commuting diagonal
self-adjoint matrices $X,Y$ such that $\|X\|\leq1$ and
\[
 \tr_d\bigl(\Ssh_{n,m}(X,Y)\bigr)
 =\tau\bigl(\Ssh_{n,m}(A,B)\bigr)\qquad(n\in N).
\]
\end{theorem}
\begin{proof}
The case $N=\varnothing$ is vacuous. Otherwise put $M=\max N$
and choose $M+1$ distinct real numbers
$x_0,\ldots,x_M\in[-1,1]$. Define
\[
 c_n=\binom{n+m}{m}^{-1}
 \tau\bigl(\Ssh_{n,m}(A,B)\bigr),\qquad 0\leq n\leq M.
\]
Each $c_n$ is real because $\Ssh_{n,m}(A,B)$ is self-adjoint: taking adjoints permutes the words in $\W_{n,m}$.
Since the nodes are distinct, the Vandermonde matrix $(x_j^n)_{0\leq n,j\leq M}$ is invertible over $\R$. Let $w_0,\ldots,w_M$ be the unique real solution of the system
\[
 \sum_{j=0}^M w_jx_j^n=c_n,\qquad 0\leq n\leq M.
\]
Set $d=M+1$ and
\[
 X=\operatorname{diag}(x_0,\ldots,x_M),\qquad
 Y=\operatorname{diag}\bigl((d w_0)^{1/m},\ldots,
                 (d w_M)^{1/m}\bigr),
\]
where the roots are real because $m$ is odd.
The matrices commute, $\|X\|\leq1$, and
$\Ssh_{n,m}(X,Y)=\binom{n+m}{m}X^nY^m$. Consequently,
\[
 \tr_d(\Ssh_{n,m}(X,Y))
 =\binom{n+m}{m}\sum_{j=0}^M w_jx_j^n
 =\tau(\Ssh_{n,m}(A,B))
\]
for every $0\leq n\leq M$.
\end{proof}

The next construction gives one matrix pair for every $n\geq0$ when $A$ has finite spectrum.

\begin{proposition}\label{thm:odd}
Let $m\geq3$ be odd, and suppose that $A$ is a self-adjoint
contraction with finite spectrum in a finite tracial von Neumann
algebra. For every self-adjoint contraction $B$, there are
self-adjoint finite matrices $X,Y$ with $\|X\|\leq1$ such that
\[
 \tr_d(\Ssh_{n,m}(X,Y))=\tau(\Ssh_{n,m}(A,B))
 \qquad\text{for every }n\geq0.
\]
\end{proposition}
\begin{proof}
Write $A=\sum_{i=1}^q a_iP_i$. For an index multiset $\alpha$ of cardinality $m$ chosen from $\{1,\ldots,q\}$, let
\begin{equation*}
 c_\alpha=\sum_{(i_1,\ldots,i_m)\sim\alpha}
 \tau(P_{i_1}BP_{i_2}B\cdots P_{i_m}B),
\end{equation*}
Here $(i_1,\ldots,i_m)\sim\alpha$ means that the tuple has multiset $\alpha$. Taking adjoints reverses the index string and conjugates the trace. Hence $c_\alpha\in\R$. As $h_n$ is symmetric, \eqref{eq:finite} becomes
\begin{equation}\label{eq:alpha}
 \tau\bigl(\Ssh_{n,m}(A,B)\bigr)
 =\frac{n+m}{m}\sum_\alpha
 c_\alpha h_n(b_{\alpha,1},\ldots,b_{\alpha,m}),
\end{equation}
where $b_{\alpha,1},\ldots,b_{\alpha,m}$ are the spectral values selected by $\alpha$, including multiplicities.

Let $R$ be the number of such multisets. For each $\alpha$, introduce a block of size $m$ and set
\begin{equation}\label{eq:cycleX}
 X_\alpha=\operatorname{diag}(b_{\alpha,1},\ldots,b_{\alpha,m}).
\end{equation}
Let $Y_\alpha$ be the real symmetric weighted adjacency matrix of the cycle
$1\! -\!2\! -\!\cdots\! -\!m\! -\!1$. Write its successive edge weights as $u_{\alpha,1},\ldots,u_{\alpha,m}$ and their product as $p_\alpha$. For $m=3$ the cycle is the ordinary triangle; it has three distinct undirected edges.

Consider a diagonal entry of
$Y_\alpha X_\alpha^{r_1}\cdots Y_\alpha X_\alpha^{r_m}$.
Expanding matrix multiplication gives a sum over length-$m$ closed walks on the cycle. Give each step an increment $+1$ or $-1$ modulo $m$. Closedness implies that the sum of the $m$ increments is divisible by $m$. Since $m$ is odd, this sum cannot equal zero. Its absolute value is at most $m$, so it must be $m$ or $-m$. Therefore every step has the same orientation. There are exactly $m$ clockwise and $m$ counterclockwise closed walks of length $m$, all of which visit every vertex once before returning. Each walk contributes the edge product $p_\alpha$. Summing over the weak compositions of $n$ distributes the powers among all vertices and gives the common scalar factor $h_n(b_{\alpha,1},\ldots,b_{\alpha,m})$. Consequently,
\begin{equation}\label{eq:cycleexact}
 \sum_{r_1+\cdots+r_m=n}
 \Tr\bigl(Y_\alpha X_\alpha^{r_1}\cdots Y_\alpha X_\alpha^{r_m}\bigr)
 =2m\,p_\alpha h_n(b_{\alpha,1},\ldots,b_{\alpha,m}).
\end{equation}

Take the direct sums $X=\bigoplus_\alpha X_\alpha$ and $Y=\bigoplus_\alpha Y_\alpha$. Their dimension is $d=mR$. By \eqref{eq:cycleexact}, the coefficient of $h_n(b_{\alpha,1},\ldots,b_{\alpha,m})$ in the normalized cyclic sum is $2mp_\alpha/d$. Set
\begin{equation*}
 p_\alpha=\frac{d}{2m}\,c_\alpha=\frac R2c_\alpha.
\end{equation*}
Choose $u_{\alpha,1}=p_\alpha$ and $u_{\alpha,j}=1$ for $2\leq j\leq m$. Then \eqref{eq:cyclic}, \eqref{eq:alpha}, and \eqref{eq:cycleexact} give
\begin{equation*}
 \tr_d\bigl(\Ssh_{n,m}(X,Y)\bigr)
 =\tau\bigl(\Ssh_{n,m}(A,B)\bigr)\qquad(n\geq0).
\end{equation*}
All blocks are real symmetric, so $X,Y$ are self-adjoint.

\end{proof}

For arbitrary output matrices, directed cycles give an exact model for all degrees when $A$ has finite spectrum.
\begin{proposition}\label{prop:directed}
Let $m\geq1$ and let $A$ have finite spectrum. There are finite matrices $X,Y$ such that $\tr_d(\Ssh_{n,m}(X,Y))=\tau(\Ssh_{n,m}(A,B))$ for every $n\geq0$. For arbitrary self-adjoint contractions $A,B$, the same construction gives simultaneous approximation for every prescribed finite set of degrees.
\end{proposition}
\begin{proof}
For finite-spectrum $A$, use the same multiset grouping \eqref{eq:alpha}. On a block of size $m$, take the diagonal matrix \eqref{eq:cycleX}, and let $Y_\alpha$ be a weighted cyclic shift: its only nonzero entries correspond to $1\to2\to\cdots\to m\to1$. There are $m$ closed directed walks of length $m$, one starting at each vertex, and all have the same edge product $p_\alpha$. Therefore the analogue of \eqref{eq:cycleexact} has $m p_\alpha$ in place of $2m p_\alpha$. With $d=mR$, choose $p_\alpha=d c_\alpha/m$ to match \eqref{eq:alpha}. For $m=1$ this is simply a one-dimensional diagonal block. For a general contraction $A$, let $A_k$ be finite-spectrum
self-adjoint contractions converging in norm to $A$. For every
fixed $n$, telescoping the occurrences of $A$ in each word gives
\[
 \left|\tau\bigl(\Ssh_{n,m}(A_k,B)-\Ssh_{n,m}(A,B)\bigr)\right|
 \leq n\binom{n+m}{m}\|A_k-A\|.
\]
Given a finite set of degrees, choose $k$ so that all these
errors are smaller than the prescribed tolerance and use the
exact directed-cycle model for $A_k$.
\end{proof}

\section{An even-order criterion}\label{sec:even}
\begin{theorem}\label{thm:even}
For a fixed even integer $m\geq2$, the following are equivalent.
\begin{enumerate}[label=\textup{(\roman*)}]
\item For every finite tracial von Neumann algebra, every pair of self-adjoint contractions, and every finite $N\subset\N_0$, the traces $\tau(\Ssh_{n,m}(A,B))$, $n\in N$, admit simultaneous contraction models.
\item For every such algebra and pair, $\Lambda_{A,B}^{(m)}(q)\geq0$ for every real polynomial $q$ nonnegative on $[-1,1]$.
\end{enumerate}
When these conditions hold, there is a positive measure $\nu_{A,B}^{(m)}$ on $[-1,1]$ satisfying
\begin{equation}\label{eq:even-meas}
 \tau\bigl(\Ssh_{n,m}(A,B)\bigr)
 =\binom{n+m}{m}\int_{-1}^1x^n\,d\nu_{A,B}^{(m)}(x),
 \qquad \nu_{A,B}^{(m)}([-1,1])=\tau(B^m)\leq1.
\end{equation}
\end{theorem}
\begin{proof}
Assume (ii). Since $m$ is even and $B$ is a self-adjoint contraction, $\tau(B^m)$ lies in $[0,1]$. Taking $f(x)=x^m/m!$ in the definition of $\Lambda$ gives
\begin{equation}\label{eq:massm}
 \Lambda_{A,B}^{(m)}(1)=\tau(B^m).
\end{equation}
Applying positivity to $\|q\|_\infty\pm q$ and using \eqref{eq:massm} gives
$|\Lambda_{A,B}^{(m)}(q)|\leq\tau(B^m)\|q\|_\infty$.
The Riesz representation theorem gives a positive measure $\nu_{A,B}^{(m)}$ on $[-1,1]$ with total mass $\tau(B^m)$. With
$f_n(x)=\frac{n!}{(n+m)!}x^{n+m}$, we obtain
\begin{align*}
 \Lambda_{A,B}^{(m)}(x^n)
 &=\frac{n!\,m!}{(n+m)!}\tau\bigl(\Ssh_{n,m}(A,B)\bigr)\\
 &=\binom{n+m}{m}^{-1}\tau\bigl(\Ssh_{n,m}(A,B)\bigr).
\end{align*}
Thus $\nu_{A,B}^{(m)}$ represents every normalized shuffle coefficient as in \eqref{eq:even-meas}. Applying Lemma \ref{lem:atomic} to the finite list $\binom{n+m}{m}x^n$ and placing the atomic masses on commuting diagonal matrices with $Y$-entries equal to one and zero, as in Corollary \ref{cor:four-measure}, gives (i).

Conversely, assume (i), and let
$q(x)=\sum_{n=0}^N c_nx^n\geq0$ on $[-1,1]$. Choosing
$f(x)=\sum_{n=0}^N c_n\frac{n!}{(n+m)!}x^{n+m}$ and expanding the coefficient of $t^m$ gives
\begin{equation}\label{eq:combLambda}
 \Lambda_{A,B}^{(m)}(q)=\sum_{n=0}^N
 \frac{c_n}{\binom{n+m}{m}}\tau\bigl(\Ssh_{n,m}(A,B)\bigr).
\end{equation}
For each positive integer $k$, apply (i) with the finite set $\{0,\ldots,N\}$ and tolerance $1/k$, obtaining a pair of self-adjoint contraction matrices $X_k,Y_k$. Since \eqref{eq:combLambda} is a finite linear combination of the specified shuffle traces,
\begin{equation}\label{eq:limitLambda}
 \Lambda_{A,B}^{(m)}(q)=\lim_{k\to\infty}
 \left.\frac{d^m}{dt^m}\tr_{d_k} f(X_k+tY_k)\right|_{t=0}.
\end{equation}
Hein\"avaara's tracial joint spectral measure theorem states that for Hermitian matrices and even $m$, a real smooth function with nonnegative $m$th derivative has nonnegative $m$th trace derivative \cite[Corollary~1.2]{HeinavaaraTJSM}. Since $q=f^{(m)}$ is nonnegative on $[-1,1]$, we localize it to apply the matrix theorem, which is stated for functions with nonnegative derivatives on $\R$. Fix $k$ and $0<\delta<1$, and replace $X_k$ temporarily by $(1-\delta)X_k$. Its spectrum lies strictly inside $(-1,1)$. Choose a smooth cutoff $\eta$ equal to one on a neighborhood of this finite spectrum, with $0\leq\eta\leq1$ and support inside $(-1,1)$. Then $\widetilde q=\eta q+(1-\eta)$ is smooth and nonnegative on the whole real line. Choose a smooth $\widetilde f$ with $\widetilde f^{(m)}=\widetilde q$ that agrees with $f$, up to a polynomial of degree less than $m$, near the spectrum of $(1-\delta)X_k$. Their derivatives through order $m$ agree in a neighborhood of the finite spectrum, so the corresponding $m$th trace derivatives at zero agree. The matrix theorem gives
\[
 \left.\frac{d^m}{dt^m}\tr_{d_k}f((1-\delta)X_k+tY_k)\right|_{t=0}\geq0.
\]
Because $f$ is a polynomial, this derivative depends continuously on the matrix $X_k$. Letting $\delta\downarrow0$ proves the same inequality with $X_k$. Thus every term on the right of \eqref{eq:limitLambda} is nonnegative. Taking the limit proves (ii).
\end{proof}

\section{The sixth-order question}\label{sec:six}
By Theorem \ref{thm:even}, the six-occurrence contraction model property is equivalent to the following assertion.
\begin{conjecture}\label{conj:six}
For self-adjoint contractions $A,B$ in every finite tracial von Neumann algebra and every real polynomial $q\geq0$ on $[-1,1]$,
\[
 \Lambda_{A,B}^{(6)}(q)\geq0.
\]
\end{conjecture}
The kernel obtained by splitting six-step cyclic paths into two three-step paths need not be positive semidefinite. For $f(x)=x^{12}$, the cyclic expansion involves the complete homogeneous polynomial $h_6$ in six spectral variables. Cutting a cyclic path of length six into two paths of length three with fixed endpoint values $a,d$ gives the scalar kernel
\begin{equation}\label{eq:K}
 K_{a,d}((b,c),(e,f))=h_6(a,b,c,d,e,f).
\end{equation}
\begin{proposition}\label{prop:indef}
The kernel \eqref{eq:K} is not positive semidefinite in general.
\end{proposition}
\begin{proof}
Choose $a=1$, $d=-1$ and the three pairs
$z_1=(1,1)$, $z_2=(1,-1)$ and $z_3=(-1,-1)$. The generating function for complete homogeneous polynomials is
\begin{equation*}
 \sum_{r\geq0}h_r(x_1,\ldots,x_6)t^r
 =\prod_{j=1}^6(1-x_jt)^{-1}.
\end{equation*}
For a six-tuple with $p$ entries equal to $1$ and $6-p$ entries equal to $-1$, the sixth complete homogeneous polynomial equals
\[
 [t^6](1-t)^{-p}(1+t)^{-(6-p)}.
\]
Evaluating this coefficient at the nine pairs $(z_i,z_j)$ gives the principal matrix
\[
 K=\begin{pmatrix}130&30&10\\30&10&30\\10&30&130\end{pmatrix}.
\]
Its determinant is
\begin{align*}
 \det K&=130(1300-900)-30(3900-300)+10(900-100)\\
 &=52000-108000+8000=-48000.
\end{align*}
Since $\det K<0$, the kernel is not positive semidefinite.
\end{proof}

\appendix
\section{Details of the fourth-order sign calculation}\label{app:signs}
For a cyclic sign pattern $\sigma_1\sigma_2\sigma_3\sigma_4$ with each
$\sigma_i\in\{+,-\}$, let
\begin{equation}\label{eq:Ssign}
 T_{\sigma_1\sigma_2\sigma_3\sigma_4}
 =\sum_{i_j\in I_{\sigma_j}}
 [\lambda_{i_1},\lambda_{i_2},\lambda_{i_3},\lambda_{i_4}]_g
 \tau(P_{i_1}BP_{i_2}BP_{i_3}BP_{i_4}B),
\end{equation}
where $g(x)=x|x|$. A cyclic rotation of the four indices does not change the
summand, and divided differences are symmetric in all their arguments.
Consequently,
\begin{equation}\label{eq:sign-group}
 T_g=4\bigl(T_{+---}+T_{++--}+T_{-+++}\bigr)
      +2T_{+-+-}.
\end{equation}
Here the all-positive and all-negative terms vanish. The coefficient four in
the first bracket counts cyclic placements of each two-sign-change pattern,
while the coefficient two counts the two alternating starting signs. The factors count marked positions of cyclic sign patterns, including repeated spectral values.

For $i\in I_+$ and $j\in I_-$, split $C_{ij}$ from \eqref{eq:C} as
\begin{align*}
 U_{ij}&=\sum_{r\in I_-}
 \frac{\lambda_i}{\lambda_i-\lambda_r}P_iBP_rBP_j,\\
 V_{ij}&=\sum_{r\in I_+}
 \frac{\lambda_j}{\lambda_j-\lambda_r}P_iBP_rBP_j.
\end{align*}
Both $U_{ij}$ and $V_{ij}$ lie in $P_i\M P_j$; thus their mixed products lie in $P_i\M P_i$. Substitution of \eqref{eq:DD1} and \eqref{eq:DD3}
into \eqref{eq:Ssign}, followed by summation over the two internal indices,
gives
\begin{align}\label{eq:pure-patterns}
 T_{+---}&=2\sum_{i\in I_+,j\in I_-}
       \frac{\tau(U_{ij}U_{ij}^*)}{\lambda_i-\lambda_j},\\
 T_{-+++}&=2\sum_{i\in I_+,j\in I_-}
       \frac{\tau(V_{ij}V_{ij}^*)}{\lambda_i-\lambda_j}.
\end{align}
Expanding the first right-hand side of \eqref{eq:pure-patterns} gives
\[
 2\sum_{i\in I_+,j\in I_-}\sum_{r,s\in I_-}
 \frac{\lambda_i^2}{(\lambda_i-\lambda_j)
 (\lambda_i-\lambda_r)(\lambda_i-\lambda_s)}
 \tau(P_iBP_rBP_jBP_sB).
\]
Cyclically rotate the four factors so that the positive index occurs first. The resulting expression is $T_{+---}$.

For the mixed pattern $++--$, fix $b_1,b_2>0$ and $a_1,a_2<0$. The coefficient of
$\tau(P_{b_1}BP_{b_2}BP_{a_1}BP_{a_2}B)$ in $T_{++--}$ is
\begin{equation}\label{eq:mixed-DD}
 [b_1,b_2,a_1,a_2]_g
 =\frac{2((b_1+b_2)a_1a_2-(a_1+a_2)b_1b_2)}
 {(b_1-a_1)(b_1-a_2)(b_2-a_1)(b_2-a_2)}.
\end{equation}
Expanding $\tau(U_{ij}V_{ij}^*)$ and $\tau(V_{ij}U_{ij}^*)$ and grouping the two cyclic orderings gives
\begin{align}\label{eq:mixed-rational}
 &\frac{2b_1a_1}{(b_1-a_1)(b_1-a_2)(a_1-b_2)}
 +\frac{2b_2a_2}{(b_2-a_2)(b_2-a_1)(a_2-b_1)}\\
 &\quad=\frac{2((b_1+b_2)a_1a_2-(a_1+a_2)b_1b_2)}
 {(b_1-a_1)(b_1-a_2)(b_2-a_1)(b_2-a_2)}.\nonumber
\end{align}
Combining the fractions in \eqref{eq:mixed-rational} gives \eqref{eq:mixed-DD}. The cyclic trace identity then yields
\begin{equation}\label{eq:mixed-pattern}
 T_{++--}=2\sum_{i\in I_+,j\in I_-}
 \frac{\tau(U_{ij}V_{ij}^*+V_{ij}U_{ij}^*)}{\lambda_i-\lambda_j}.
\end{equation}
For repeated nodes within one sign class, take the continuous limits of the divided differences in \eqref{eq:DD1}--\eqref{eq:DD3}.
Adding \eqref{eq:pure-patterns} and \eqref{eq:mixed-pattern}, we obtain
\begin{align*}
 T_{+---}+T_{++--}+T_{-+++}
 &=2\sum_{i\in I_+,j\in I_-}\frac{1}{\lambda_i-\lambda_j}
 \tau\bigl((U_{ij}+V_{ij})(U_{ij}+V_{ij})^*\bigr)\\
 &=Q_{\rm adj}.
\end{align*}
Its contribution to \eqref{eq:sign-group} is $4Q_{\rm adj}$.

For the alternating pattern, let $a_1,a_2<0<b_1,b_2$.
Writing \eqref{eq:DD2} with a common positive denominator gives
\begin{equation}\label{eq:alt-scalar}
 [a_1,a_2,b_1,b_2]_g
 =2\frac{(-a_1-a_2)b_1b_2+(b_1+b_2)a_1a_2}
 {(b_1-a_1)(b_1-a_2)(b_2-a_1)(b_2-a_2)}.
\end{equation}
The two terms in the numerator of \eqref{eq:alt-scalar} give, respectively,
\begin{align*}
&2\sum_{\substack{i,j\in I_-\\r,s\in I_+}}
\frac{(-\lambda_i-\lambda_j)\lambda_r\lambda_s}
{(\lambda_i-\lambda_r)(\lambda_j-\lambda_r)
 (\lambda_i-\lambda_s)(\lambda_j-\lambda_s)}
\tau(P_iBP_rBP_jBP_sB)\\
&\qquad=2\sum_{i,j\in I_-}(-\lambda_i-\lambda_j)
   \tau\bigl(D^-_{ij}(D^-_{ij})^*\bigr),
\end{align*}
\begin{align*}
&2\sum_{\substack{i,j\in I_+\\r,s\in I_-}}
\frac{(\lambda_i+\lambda_j)\lambda_r\lambda_s}
{(\lambda_i-\lambda_r)(\lambda_j-\lambda_r)
 (\lambda_i-\lambda_s)(\lambda_j-\lambda_s)}
\tau(P_iBP_rBP_jBP_sB)\\
&\qquad=2\sum_{i,j\in I_+}(\lambda_i+\lambda_j)
   \tau\bigl(D^+_{ij}(D^+_{ij})^*\bigr).
\end{align*}
The sum of these expressions is $2Q_{\rm alt}$.
\begin{equation}\label{eq:alt-pattern}
 T_{+-+-}=2Q_{\rm alt}.
\end{equation}
Equations \eqref{eq:sign-group}, \eqref{eq:alt-pattern}, and the
preceding two-change calculation yield
$T_g=4(Q_{\rm adj}+Q_{\rm alt})$, as asserted in \eqref{eq:Gram}.

\bigskip
 
\end{document}